\documentclass[11pt]{article}
\usepackage{amsmath, amsthm, amssymb}
\begin{document}
\title{Asymptotics of the Gaussian Curvatures of the 
Canonical Metric on the Surface}

\author{Zheng Huang}
\date{}
\newtheorem{theorem}{Theorem}[section]

\newtheorem{pro}[theorem]{Proposition}
\newtheorem{cor}[theorem]{Corollary}
\newtheorem{lem}[theorem]{Lemma}
\newtheorem{rem}[theorem]{Remark}

\newcommand{\WP}{Weil-Petersson}
\newcommand{\TS}{Teichm\"{u}ller space}
\newcommand{\ms}{moduli space}
\newcommand{\cs}{conformal structure}
\newcommand{\Bm}{Bergman metric}
\newcommand{\cm}{canonical metric}
\newcommand{\hq}{holomorphic quadratic}
\newcommand{\RS}{Riemann surface}
\newcommand{\Sc}{sectional curvature}
\newcommand{\cd}{complex dimension}
\newcommand{\Bd}{Beltrami differential}
\newcommand{\ts}{tangent space}
\newcommand{\Gc}{Gaussian curvature}
\newcommand{\Ad}{Abelian differential}
\newcommand{\pc}{pinching coordinate}
\newcommand{\mcg}{mapping class group}

\maketitle
\begin{abstract}
We study the canonical metric on a compact Riemann surface of genus at least two. This natural 
metric is the pullback, via the period map, from the Euclidean metric on the Jacobian variety of 
the surface. While it is known that the canonical metric is of nonpositive curvature, we show that 
its Gaussian curvatures are not bounded away from zero nor negative infinity when the surface is 
close to the compactification divisor of Riemann's moduli space. 
\end{abstract}
\tableofcontents
\footnotetext{Zheng Huang: Department of Mathematics, 
University of Michigan, Ann Arbor, MI 48109, USA. Email: zhengh@umich.edu}
\footnotetext{AMS subject classification 30F60, 32G15}
\section {Introduction}
Throughout this paper $\Sigma$ is a compact, smooth, oriented, closed {\RS} of 
genus $g \ge 1$. The period map $p: \Sigma \rightarrow J_{\Sigma}$ embeds 
the surface $\Sigma$ to its Jacobian $J_{\Sigma}$. The pullback metric of the 
Euclidean metric on $J_{\Sigma}$ via this period map thus defines the so-called 
{\it {\cm}} or {\it {\Bm}} on $\Sigma$, denoted by $\rho_B$. This metric 
$\rho_B$ is of nonpositive {\Gc}, and when $g \ge 2$, the curvature vanishes if 
and only if the surface is hyperelliptic and only at $2g+2$ Weierstrass points 
(\cite {GR1}, \cite {L}), in other words, amazingly, the {\Gc}s characterize 
hyperelliptic surfaces.

In this short article, we study the asymptotics of the {\Gc}s of the {\cm}, when the 
surface is degenerating towards the compactification divisor of the {\ms}. We show 
that, using {\pc}s and classical techniques on asymptotics of {\Ad}s,
\begin{theorem}
If $g(\Sigma) > 2$, the {\Gc}s of the {\cm} are not bounded from below, nor from above, 
by any negative constants independent of the conformal structure of the surface.
\end{theorem}

This unboundedness result is not surprising, however, it is not clear about the rates of the asymptotic 
behavior of the curvatures.  We will make the theorem more precise by providing quantitative 
estimates on the {\Gc}s in theorems 2.4 and 2.6.

This paper is organized as follows. We give a quick exposition of {\ms} and {\pc} in 
$\S 2.1$, and discuss the {\cm} on a compact {\RS} in $\S 2.2$, then proceed 
to calculate the asymptotics of its {\Gc} and prove theorem 1.1 in $\S 2.3$. 

The author acknowledges several enjoyable and helpful conversations 
with Jiaqiang Mei and Ralf Spatzier.

\section{Proof of the Main Theorem}
We give some background on {\ms}, {\pc}s, Masur's frame of regular quadratic 
differentials and its modification in $\S 2.1$; we introduce the {\cm} on a 
{\RS} in $\S 2.2$; and then calculate the asymptotics of its {\Gc}s in 
$\S 2.3$, where we prove theorem 1.1 as a consequence of our calculation.

\subsection{Differentials and Pinching Coordinates} 
We recall that {\TS} ${\mathcal {T}}_{g}$ is the space of {\cs}s on 
$\Sigma$, where two {\cs}s $\sigma$ and $\rho$ are equivalent if there is a 
biholomorphic map between $(\Sigma,\sigma)$ and $(\Sigma,\rho)$ in the 
homotopy class of the identity. Riemann's {\ms} $\mathcal{M}_g$ of {\RS}s is 
obtained as the quotient of {\TS} by the {\mcg}. {\TS} is a complex manifold of 
complex dimension $3g - 3 > 1$ when $g>1$,  and the co{\ts} at $\Sigma$ is 
identified with $QD(\Sigma)$, the space of {\hq} differentials while the tangent 
space consists of harmonic {\Bd}s.

We will use the {\pc}s to describe the degeneration of the surface, following 
{\cite{EM}, \cite{Ma}}, and the construction of Wolpert (\cite {Wp03}). 
To this goal, we consider the Deligne-Mumford compactification of $\mathcal{M}_g$, 
and any element of the compactifying divisor is a {\RS} with nodes (\cite{B}).

Let $\Sigma_0$ be a {\RS} obtained from pinching disjoint, nonhomotopic and 
noncontractible simple closed curves $\gamma_1, \cdots, \gamma_m$ on 
$\Sigma$ to a point, then $\Sigma_0$ has paired punctures 
$a_1, b_1, \cdots, a_m, b_m$ corresponding to the curves 
$\gamma_1, \cdots, \gamma_m$. Thus $\Sigma_{0}$ corresponds to a point at the 
compactifying divisor. For simplicity, we assume $\Sigma_0$ does not have spheres 
or tori as components.

Let $a_j$ and $b_j$ be two paired punctures of $\Sigma_0$, then one node on 
$\Sigma_0$ is obtained by identifying $a_j$ and $b_j$. Let 
$D = \{ z \in C: |z| < 1 \}$, and $U_j$ and $V_j$ be small disjoint 
neighborhoods of $a_j$ and $b_j$, respectively, where $j=1,\cdots,m$, we have 
local coordinates $z_j: U_j \rightarrow D$ and $w_j: V_j \rightarrow D$ such 
that $z_j(a_j)=0$ and $w_j(b_j)=0$. Given $t_j \in C$ with $0 < |t_j| < 1$, we 
obtain a new {\RS} with nodes by removing the discs $\{ |z_j| < |t_j| \}$ and 
$\{ |w_j| < |t_j| \}$ and then identifying $z_j$ with $w_j = {t_j}/{z_j}$. This is 
the process of opening the node. Therefore 
$\{ t=(t_1, \cdots, t_m) \in C^m : |t_j| < 1 \}$ forms the complex parameters for 
opening the nodes.

We choose a basis of {\Bd}s $\{ \mu_1, \cdots, \mu_{3g-3-m} \}$ for the space 
$T_{\Sigma_0}{\mathcal{T}}(\Sigma_{0})$ with support in 
$\Sigma_{0} \backslash \bigcup_{j=1}^{m}(U_j \cup V_j)$, then 
$\tau = (\tau_1, \cdots, \tau_{3g-3-m})$ in an open set $U \subset C^{3g-3-m}$ 
about the origin provides local coordinates for the neighborhood of $\Sigma_0$ 
in ${\mathcal{T}}(\Sigma_{0})$. These parameters $(t,\tau)$ define a local covering for the 
compactified {\ms} $\bar{\mathcal{M}}_g$ around $\Sigma_0$, and we call them 
{\it {\pc}s}. For $\{ |t_j| > 0 \}$, we have a compact {\RS} 
$\Sigma_{(t,\tau)}$. In these coordinates, the compactifying divisor 
$D = \bar{\mathcal{M}}_g - {\mathcal{M}}_g$ is represented by the equation $\prod_{j=1}^{m}t_j=0$.

For $(t, \tau)$ near $(0,0)$, let $\gamma_1, \cdots, \gamma_m$ be $m$ short 
geodesics on the surface. Among those curves, we assume 
$\gamma_1, \cdots, \gamma_{m_1}$ are nonseparating ($m_1 \le m$) while the 
rest are separating curves. Masur constructed a local frame 
$\{ \phi_i \}_{i=1}^{3g-3}(t,\tau)$ of regular quadratic differentials which 
is holomorphic in the {\pc}s $(t,\tau)$ and represents the canonical coframe 
$(dt,d\tau)$ in {\pc}s. In other words, 
$\{ \phi_{j}dz^2, \phi_{\nu}dz^2 \}_{1 \le j \le m; \nu \ge m+1}$ is a basis 
dual to the basis formed by 
$\{ \partial/\partial t_j, \partial/\partial \tau_{\nu} 
\}_{1 \le j \le m; \nu \ge m+1}$ {\cite{Ma}}.

More specifically, by {\cite{Ma}}, on $U_j$, for $1 \le i \le m$ and $\nu \ge m+1$, 
\begin{eqnarray}
\phi_{i}(z_j,t,\tau) = -{\frac{t_i}{\pi}}\{ {\frac{\delta_{ij}}{z_{j}^2}} 
+ a_{-1}(z_j,t,\tau) + {\frac{1}{z_{j}^2}}\sum_{k=1}^{\infty}
[\frac{t_j}{z_j}]^{k}t_{j}^{m(k)}a_{k}(t,\tau)\} 
\end{eqnarray}
where $m(k) \ge 0$, and $a_{-1}$ with at most a simple pole at $z_{j}=0$ 
and $a_k$ is holomorphic in $t$ and $\tau$ for $k \ge 1$. And
\begin{eqnarray}
\phi_{\nu}(z_j,t,\tau) = \phi_{\nu}(z_j,0,0) + 
{\frac{1}{z_{j}^2}}\sum_{k=1}^{\infty}
[\frac{t_j}{z_j}]^{k}t_{j}^{m'(k)}b_{k} + 
\sum_{k=-1}^{\infty}z_{j}^{k}c_{k}
\end{eqnarray}
where $m'(k) \ge 0$, and $\phi_{\nu}(z_j,0,0)$ has at most a simple pole, 
functions $b_k = b_k(t,\tau)$ and $c_k = c_k(t,\tau)$ are holomorphic 
around 0. We have similar formulas for regular quadratic differentials on 
$V_j$ in $(w_j,t,\tau)$ coordinates. Therefore, for $1 \le i \not= j \le m$, 
and if $K$ is any compact set on the surface that contains no singularity, 
we find that for $1 \le j \le m$:
\begin{eqnarray}
|\phi_{j}| &\sim& |t_j|/|z_j|^2 \mbox{ on } U_j \nonumber \\
|\phi_{j}| &=& O(|t_j|/|z_i| + |t_{i}t_{j}|/|z_i|^3) \mbox{ on } U_i, 
1 \le i \not= j \le m \nonumber \\
|\phi_{j}| &=& O(|t_j|) \mbox{ on } K 
\end{eqnarray}
where $A \sim B$ means $A/C \le B \le CA$ for some constant $C > 0$.

Habermann and Jost (\cite{HJ}) modified Masur's local frame to obtain a new coframe 
$\{ \psi_{j}dz^2, \psi_{\nu}dz^2 \}_{1 \le j \le m; \nu \ge m+1}$:
\begin{center}
$ \psi_{j}(t,\tau) = \left \{
\begin{array}{cc}
\phi_{j}(t,\tau) - t_{j}\sum_{k=1}^{3m}\lambda_{jk}(t,\tau)
\phi_{k}(t,\tau) & \mbox{ if } 1 \le j \le m \\
\sum_{k=m+1}^{3g-3}\lambda_{jk}(t,\tau)\phi_{k}(t,\tau) , & \mbox{ if } 
j>m
\end{array} \right. $
\end{center}
where functions $\lambda_{jk}(t,\tau)$ are holomorphic in $(t, \tau)$ while 
$\lambda_{jk}(0,0) = \delta_{jk}$ if $j>m$. We notice that each node corresponds 
to three quadratic differentials: one differential corresponds to the 
$t$-direction while two other correspond to change of the position of the paired 
punctures associated to this node. We collect the properties of this modified frame from 
(\cite{HJ}), for $1 \le j \le m$:
\begin{eqnarray}
|\psi_{j}| &\sim& |t_j|/|z_j|^2 \mbox{ on } U_j \nonumber \\
|\psi_{j}| &=& O(|t_j|) \mbox{ on } U_i, 1 \le i \not= j \le m \nonumber \\
|\psi_{m+2j-1}| &\sim& |\psi_{m+2j}| \sim (1/|z_j|) \mbox{ on } U_j 
\nonumber \\
|\psi_{\nu}| &=& O(|t_j|) \mbox{ on } U_j, \mbox{ for } 3m+1 \le \nu 
\nonumber \\
|\psi_{j}| &=& O(|t_j|) \mbox{ on } K 
\end{eqnarray}
\subsection {The Canonical Metric}
We now introduce the {\cm} on a {\RS}, and investigate the asymptotic behavior 
of its {\Gc}s.

On a compact {\RS} $\Sigma$ of genus $g>1$, the dimension of the space of 
{\Ad}s of the first kind, or holomorphic one forms, is $g$. There is a natural 
pairing of {\Ad}s defined on this space:
\begin{center}
$<\mu, \nu> = {\frac{\sqrt{-1}}{2}}\int_{\Sigma}\mu\wedge\bar{\nu}$
\end{center}
Let $\{\omega_1, \omega_2, \cdots, \omega_g \}$ be a basis of {\Ad}s, normalized 
with respect to the $A$-cycles of some symplectic homology basis 
$\{A_i, B_i \}_{1 \le i \le g}$, i.e., $\int_{A_i}\omega_j = \delta_{ij}$. 
Thus the period matrix $\Omega_{ij} = \int_{B_i}\omega_j$. One finds that, since 
not all {\Ad}s vanish at the same point according to Riemann-Roch, the period 
matrix is then symmetric with positive definite imaginary part:
$Im \Omega_{ij} = <\omega_i, \omega_j>$ (\cite{FK}) .

The {\cm} $\rho_B$ on surface $\Sigma$ is the metric associated to the $(1,1)$ 
form given by 
\begin{center}
$ {\frac{\sqrt{-1}}{2}}\sum_{i,j=1}^{g}
(Im\Omega)_{ij}^{-1}\omega_{i}(z){\bar{\omega}}_{j}(\bar{z})$.
\end{center}

It is not hard to see that this metric is the pull-back of the Euclidean metric 
from the Jacobian variety $J(\Sigma)$ via the period map (\cite {FK}). 
\begin{rem}
It is easy to see that the area of the surface $\Sigma$ with respect to the {\cm} 
is a constant, i.e., $\int_{\Sigma}\rho_B = g$. Sometimes the {\cm} is also 
refered to ${\frac{\rho_B}{g}}$ to unify the surface area. 
\end{rem}

Grauert-Reckziegel (\cite {GR1}), and Lewittes (\cite {L}) showed the {\Gc}s 
are nonpositive, when $g \ge 2$, and Lewittes determined further that $K_c(p) = 0$ 
for some $p \in \Sigma$ if and only if $\Sigma$ is hyperelliptic and $p$ is one of 
the $2g+2$ classical Weierstrass points of $\Sigma$. One can perturb this metric to 
obtain a metric of negative curvature on a compact {\RS} in an elementary fashion, 
without using the uniformization theorem (see \cite {GR2}).

As in $\S 2.1$, let $(t,\tau)=(t_1,\cdots,t_m,\tau_1,\cdots,\tau_{3g-3-m})$ be 
the {\pc}s, and $\Sigma_0$ be the surface corresponding to pinching $m$ short 
curves $\gamma_1,\cdots,\gamma_m$ on $\Sigma$, moreover, the curves 
$\gamma_1,\cdots,\gamma_{m_1}$ are nonseparating with $1 \le m_1 \le m$. Now 
we have paired punctures $a_1,b_1,\cdots,a_m,b_m$.

We adapt some notations from {\cite{HJ}}, defining sets: 
\begin{center}
$B_{j,t_j} = \{ z_j: |t_j|^{1/2} < |z_j| < 1 \} $\\
${\tilde{B}}_{j,t_j} = \{ w_j: |t_j|^{1/2} < |w_j| < 1 \} $
\end{center}
Let $\rho_{B}^0$ be the {\cm} on $\Sigma_0$ and $\rho_{B}(t,\tau)$ be the {\cm} 
on the surface $\Sigma_{(t,\tau)}$, we define $\rho_{j}(z_j,t,\tau)$ by 
\begin{center}
$\rho_{B}(t,\tau) = \rho_{j}(z_j,t,\tau)dz_{j}d{\bar{z}}_j$ on $B_{j,t_j}$.  
\end{center}
\begin{pro}(\cite{HJ})
\begin{itemize}
\item
On any compact set 
$K \subset \Sigma_0 \backslash \{ a_1,b_1,\cdots,a_m,b_m \}$, 
\begin{center}
$\rho_{B}(t,\tau) \rightarrow \rho_{B}^0$ as $(t,\tau) \rightarrow (0,0)$;
\end{center}
\item
For $1 \le j \le m_1$, there exist positive constants $c_j$ such that when  
$z_j \in B'_{j,t_j} = \{  z_j: |t_j|^{1/2} < |z_j| < 
|log|t_j||^{-1/2} \} \subset B_{j,t_j}$, 
\begin{eqnarray}
{\frac{1}{c_{j}|log|t_j|||z_j|^2}} \le \rho_{j}(z_j,t,\tau) 
\le {\frac{c_j}{|log|t_j|||z_j|^2}};
\end{eqnarray} 
and when $z_j \in B''_{j,t_j} = \{  z_j: |log|t_j||^{-1/2} < |z_j| < 1 \} 
\subset B_{j,t_j}$, 
\begin{eqnarray}
{\frac{1}{c_j}} \le \rho_{j}(z_j,t,\tau) \le {c_j};
\end{eqnarray} 
\item
For $z_j \in B_{j,t_j}$ and $m_1 \le j \le m$, we have
\begin{eqnarray}
{\frac{1}{c_j}} \le \rho_{j}(z_j,t,\tau) \le {c_j}.
\end{eqnarray} 
\end{itemize}
\end{pro}

We note here the pinching region corresponding to a nonseparating node becomes 
long and thin, so the above proposition actually compares the {\cm} to a flat 
long and thin cylinder for $z_j \in B'_{j,t_j}$.
\begin{rem}
There is a so-called Arakelov metric on a {\RS} whose {\Gc} is proportional to 
the {\cm}, studied by Jorgenson (\cite {Jg}), Wentworth (\cite {We}), and others. 
\end{rem}

\subsection{Asymptotics of the Gaussian Curvatures}
We now calculate the {\Gc}s $K_c$ of the {\cm} on degenerating surface $\Sigma$ 
where various simple closed curves are shortening. For the sake of exposition, 
at this point, we assume {\it the genus $g > 2$ and there is exactly a nonseparating 
simple closed curve $\gamma_1$ is shortening}.

Let $(t,\tau) = (t_1, \tau_1, \cdots, \tau_{3g-4})$ be the {\pc}s in this 
situation. We note that different $t$ means different {\cs}. To consider 
asmptotics of the {\Gc}s $K_c$, we fix $t_1$, hence fix {\cs} on $\Sigma$, then 
study how $K_c$ depends on $t_1$.

The {\Gc} is given by
\begin{center}
$K_c (z,t_1,\tau) = -{\frac{2}{\rho_{B}(z,t_1,\tau)}}
{\frac{\partial^2}{\partial z\partial \bar{z}}}log\rho_{B}(z,t_1,\tau)$.
\end{center}
\begin{theorem}
Let $g > 2$, and $z \in B''_{1,t_1} = \{ z: |log|t_1||^{-1/2} < |z| < 1 \}$, 
then $K_c(z) \sim log|t_1|$ as $z \rightarrow |log|t_1||^{-1/2}$ in $B''_{1,t_1}$. 
Here
$K_c(z) \sim log|t_1|$ means $|{\frac{K_c(z)}{log|t_1|}}|$ is comparable to 
a positive constant independent of $t_1$, for $|t_1|$ small. Consequently, the 
{\Gc} of the {\cm} on $\Sigma$ is unbounded.
\end{theorem}
\begin{proof}
We assume meromorphic 1-forms 
$\omega_1(t_1,\tau), \cdots, \omega_g(t_1,\tau)$ form a basis for the regular 
one forms on surface $\Sigma_{(t_1,\tau)}$, for $(t_1,\tau)$ near $(0,0)$.

Let $A=(a^{ij})_{1 \le i, j \le g}$ be the inverse matrix of 
$Im\Omega(t_1,\tau)$, where $\Omega$ is the period matrix for the basis 
$\omega_1(t_1,\tau), \cdots, \omega_g(t_1,\tau)$, then the {\cm} can be written 
as $\rho_B (z, t_1, \tau) = \sum_{i,j=1}^{g}a^{ij}\omega_i(z)
{\bar{\omega}}_j(\bar{z})$. The asymptotics of the period matrix has been 
extensively studied by Fay (\cite {Fy}), Yamada (\cite {Ya}), and others. 
One finds that (\cite {Fy} or \cite {Ya}), when the surface is developping a 
single, nonseparating node,
\begin{eqnarray*}
<\omega_1, \omega_1> &=& c_{-1}|log|t_1|| + c_0 +c_1|t_1| + O(|t_1|^2); \\
<\omega_j, \omega_j> &\sim& 1 + O(|t_1|), j \not= 1; \\
|<\omega_1, \omega_j>| &=& O(1), j \not= 1,
\end{eqnarray*}
where $c_{-1} \not= 0, c_0, c_1$ are constants independent of $t_1$.

One writes
\begin{center}
$\omega_i(t_1,\tau) = \omega_{i1}(z,t_1,\tau)dz_1$ on 
$A_{1,t_1} = \{z: |t_1| < |z| < 1 \}$, for $1 \le i \le g$,
\end{center}
where 
\begin{eqnarray}
\omega_{i1}(z,t_1,\tau) = \sum_{k=-1}^{\infty}b_{i1,k}(t_1,\tau)z^k + 
\sum_{k=-\infty}^{-2}b_{i1,k}(t_1,\tau){\frac{1}{t_1}}({\frac{z}{t_1}})^k
\end{eqnarray}
with $b_{11,-1}(0,\tau) = 1$, and $b_{i1,-1}(0,\tau)= 0$ for $1 < i \le g$, 
and $b_{i1,k}(t_1,\tau)$ are holomorphic in both $t_1$ and $\tau$.

A straightforward calculation shows that the {\Gc} is 
\begin{eqnarray}
K_c(z,t_1,\tau) &=& {\frac{-2}{\rho_B^3}}\{(\sum_{i,j=1}^{g}a^{ij}\omega_i(z)
{\bar{\omega}}_j(\bar{z}))(\sum_{ij=1}^{g}a^{ij}
{\frac{\partial \omega_{i}(z)}{\partial z}}
{\frac{\partial \bar{\omega}_{j}(\bar{z})}{\partial \bar{z}}})\nonumber \\
&-& (\sum_{ij=1}^{g}a^{ij}{\frac{\partial \omega_{i}(z)}{\partial z}}
\bar{\omega}_{j}(\bar{z}))(\sum_{ij=1}^{g}a^{ij}\omega_{i}(z)
{\frac{\partial \bar{\omega}_{j}(\bar{z})}{\partial \bar{z}}})\}
\end{eqnarray}

One notices that $\rho_B$ is bounded independent of $t_1$ in $B''_{1,t_1}$, as 
indicated in (6). Also one recalls that the matrix $Im\Omega$ is symmetric, 
positive definite, hence its inverse matrix $A=(a^{ij})$ satisfies: 
\begin{center}
$a^{11} > 0$, and $a^{11} \sim |log|t_1||^{-1}$, \\
$a^{jj} > 0$ and $a^{jj} = O(1)$ for $j \not=1$, \\
$a^{jk} = O(|log|t_1||^{-1})$ for $j \not= k$. 
\end{center} 
Therefore, using these estimates, we have
\begin{eqnarray*}
K_c(z,t_1,\tau) &=& {\frac{-2}{\rho_B^3}}\{\sum_{j=2}^{g}a^{11}a^{jj}
|\omega_1 \omega_{j,z}-
\omega_{1,z}\omega_j|^2 \\
&+& \sum_{j=2}^{g}(a^{1j})^2 |\omega_1 \omega_{j,z}-
\omega_{1,z}\omega_j|^2 + 
\mbox{ lower order terms}\} \\
&\sim& \sum_{j=2}^{g}a^{11}a^{jj}|\omega_1 \omega_{j,z}-
\omega_{1,z}\omega_j|^2 + \mbox{ lower order terms}
\end{eqnarray*}
Now using the expansion (8) and above estimates on $A=(a^{ij})$, one finds that, 
in this region $B''_{1,t_1}$, for $j \not= 1$,
\begin{eqnarray*}
|K_c(z,t_1,\tau)| &\sim& {\frac{1}{|log|t_1||}}
({\frac{|b_{j1,0}|^2 |b_{11,-1}|^2}{|z|^4}}) +O(1). 
\end{eqnarray*}
Here $b_{j1,0}(t_1,\tau)$ is the constant term in the expansion (8) of 
$\omega_{j1}$ where $j \not= 1, k = 0$, so $|b_{j1,0}|$ is bounded (not zero), 
independent of $t_1$, since $\omega_j(t_1,\tau)$ is close to a basis element on 
the space of regular one forms on a surface of genus $g-1$ when $t_1$ is 
sufficiently small, for $j \not= 1$ (see \cite {Fy} or \cite {Ya}). And 
$|b_{11,-1}(t_1,\tau)|$ is bounded (positive), indendent of $t_1$, as 
$b_{11,-1}(0,\tau) = 1$ and $b_{11,-1}(t_1,\tau)$ is holomorphic in $t_1$ and 
$\tau$. Therefore, 
\begin{eqnarray*}
|K_c(z,t_1,\tau)| &\sim& {\frac{1}{|log|t_1||}}({\frac{1}{|z|^4}})+O(1) 
\end{eqnarray*}
where $z \in B''_{1,t_1} = \{ z: |log|t_1||^{-1/2} < |z| < 1 \}$. Thus, we have 
\begin{itemize}
\item
$|K_c(z,t_1,\tau)|$ is bounded when $z$ is close to the outer circle $|z|=1$; 
\item
$|K_c(z,t_1,\tau)| \sim |log|t_1||$, if $z$ is close to the inner circle 
$|z| = |log|t_1||^{-1/2}$.
\end{itemize}
This completes the proof.
\end{proof}
\begin{rem}
We assume that $g \ge 3$ in theorem 2.4. If $g=2$, then the shrinking of a 
nonseparating closed curve on $\Sigma$ will reduce the surface to a torus, 
where the space of {\Ad}s is one dimensional. 
\end{rem}
Theorem 2.4 proves the first half of theorem 1.1, i.e., we proved the {\Gc}s 
of the {\cm} on $\Sigma$ are unbounded by estimating the {\Gc}s in the outer 
annuli region $B''_{1,t_1} = \{ z: |log|t_1||^{-1/2} < |z| < 1 \}$.

Now we continue to calculate $K_c(z,t_1,\tau)$ in the inner annuli region 
$B'_{1,t_1} = \{ z: \sqrt{|t_1|} < |z| < |log|t_1||^{-1/2} \}$, where $z$ is 
closer to the node. The key difference of the estimates in $B'_{1,t_1}$ and 
$B''_{1,t_1}$ is that the {\cm} $\rho_B$ is no longer bounded in $B'_{1,t_1}$. 
Now the inner circle $|z| = |log|t_1||^{-1/2}$ of $B''_{1,t_1}$ becomes the 
outer circle of $B'_{1,t_1}$. The inner circle of $B'_{1,t_1}$ is now 
$|z| = \sqrt{|t_1|}$.   
\begin{theorem}
Let $g > 2$, and $z \in B'_{1,t_1} = \{ z: \sqrt{|t_1|} < |z| < 
|log|t_1||^{-1/2} \}$, then $K_c(z) \sim log|t_1|$ as $z \rightarrow 
|log|t_1||^{-1/2}$, and $|K_c| = O(|t_1|log^{2}|t_1|)$ as $z \rightarrow 
\sqrt{|t_1|}$ in $B'_{1,t_1}$. 
\end{theorem}
\begin{proof}
From formula (5), in $B'_{1,t_1}$, the {\cm} 
$\rho_B \sim {\frac{1}{|z|^2 |log|t_1||}}$, hence we have 
${\frac{1}{\rho_B^3}} \sim |z|^6 |log|t_1||^3$.
Applying this to 
formula (9) and following the calculation in proving theorem 2.4, we find that
\begin{eqnarray*}
|K_c(z,t_1,\tau)| &\sim&  |z|^6 |log|t_1||^3 \{\sum_{j=2}^{g}a^{11}a^{jj}
|\omega_1 \omega_{j,z}-
\omega_{1,z}\omega_j|^2 \\
&+& \sum_{j=2}^{g}(a^{1j})^2 |\omega_1 \omega_{j,z}-
\omega_{1,z}\omega_j|^2 + 
\mbox{ lower order terms}\} \\
&\sim& |z|^6 |log|t_1||^3 \{{\frac{1}{|log|t_1||}}{\frac{1}{|z|^4}}(1+o(1))+O(1)\} 
\end{eqnarray*}
where $z \in B'_{1,t_1} = \{ z: \sqrt{|t_1|} < |z| < |log|t_1||^{-1/2} \}$. 
Thus, 
\begin{itemize}
\item
$|K_c(z,t_1,\tau)| \sim |log|t_1||$ if $z$ is close to the outer circle 
$|z| = |log|t_1||^{-1/2}$;
\item
$|K_c(z,t_1,\tau)| = O(|t_1|log^{2}|t_1|)$ when $z$ is close to inner 
circle $|z| = \sqrt{|t_1|}$.
\end{itemize}
\end{proof}

Theorems 2.4 and 2.6 immediately imply Theorem 1.1, where we claimed that the 
{\Gc}s of the {\cm} are not bounded from below, nor from above by any negative 
constants, if one deforms the conformal structure of the surface.
\begin{rem}
It is not hard to see that $|K_c(z,t_1,\tau)|$ cannot blow up at the rate of 
$|log|t_1||$ uniformly in region $B'_{1,t_1}$, from Gauss-Bonnet. In fact, a 
straightforward calculation shows
\begin{eqnarray*}
|\int_{B'_{1,t_1}}log{\frac{1}{|t_1|}}\rho_B| &\sim& 
|\int_{B'_{1,t_1}}{\frac{|log|t_1||}{|z|^2 |log|t_1||}}|dz|^2| \\
&\sim& |log|t_1|| - log|log|t_1||.
\end{eqnarray*} 
\end{rem}

We end this section with a brief discussion of the case where a single separating 
closed curve on $\Sigma$ is shortening. A such node will then disconnect the 
surface into two components $\Sigma_1$ and $\Sigma_2$, with 
$g(\Sigma_1) + g(\Sigma_2) = g(\Sigma)$. Since we are still assuming none of the 
limiting components are torus or sphere, we actually assume $g(\Sigma) \ge 4$.

In this case, one notices that, from formula (7), the metric $\rho_B$ stays 
bounded. Meanwhile, the upper $g(\Sigma_1) \times g(\Sigma_1)$ block of the matrix 
$Im\Omega$ is approaching $Im\Omega_1$, and the lower 
$g(\Sigma_2) \times g(\Sigma_2)$ block of the matrix $Im\Omega$ is approaching 
$Im\Omega_2$, where $\Omega_1$ and $\Omega_2$ are period matrices of the surfaces 
$\Sigma_1$ and $\Sigma_2$, respectively, for small $(t,\tau)$.

As we follow the proof of theorem 2.4, for the basis of meromorphic abelian 
differentials $\{\omega_1, \cdots, \omega_g\}$ and their expansions in (8), we 
find that $b_{i1,-1}(0,\tau) = 0$, and $b_{i1,-1}(t_1,\tau)$ are holomorphic in 
$t_1$ and $\tau$, for all $1 \le i \le g$. Then one sees that the 
{\Gc}s stay bounded, independent of $t_1$.

\begin{rem}
The {\cm} on the surface naturally induces an Hermitian metric of the {\WP} type on {\TS}. It is 
incomplete (\cite{HJ}), and we shall re-estabilish this metric via an approach of calculus of 
varaition in a separate paper. 
\end{rem}


\begin{thebibliography}{99}
\bibitem {B} L. Bers, {Spaces of degenerating {\RS}s}, Discontinous 
Groups and Riemann Surfaces, {\em Ann. of Math Studies}, {\bf 79} Princeton 
University Press, Princeton, New Jersey 1974

\bibitem {EM} C. Earle, A. Marden, {Geometric complex coordinates for 
{\TS}}, {\em Manuscript} 

\bibitem {FK} H. Farkas, I. Kra, {{\RS}s}, {\em Graduate Texts in 
Mathematics}, {\bf 71}, New York, Springer-Verlag, 1980

\bibitem {Fy} J. Fay, {Theta functions on {\RS}s}, {\em Lecture Notes in 
Mathematics}, {\bf 352}, New York, Springer-Verlag, 1973

\bibitem {GR1} H. Grauert, H. Reckziegel, {Hermitesche Metriken und normale 
Familien holomorpher Abbildungen}, {\em Math. Z.} {\bf 89} (1965) 108-125

\bibitem {GR2} H. Grauert, R. Remmert, {Plurisubharmonische funktionen in 
komplexen r\"{a}ume}, {\em Math. Z.} {\bf 65} (1956) 175-194

\bibitem {HJ} L. Habermann, J. Jost, {Riemannian metrics on {\TS}}, 
{\em Manusciripta Math.} {\bf 89} (1996) No. 3, 281-306

\bibitem {Jg} J. Jorgenson, {Asymptotic behavior of Falting's delta function}, 
{\em Duke Math. J.}, {\bf 61} (1990) 221-254 

\bibitem {L} J. Lewittes, {Differentials and metrics on {\RS}s}, 
{\em Trans. Amer. Math. Soc.} {\bf 139} (1969) 311-318

\bibitem {Ma} H. Masur, {The extension of the {\WP} metric to the 
boundary of {\TS}}, {\em Duke Math J}. {\bf 43} (1976) 623-635

\bibitem {We} R. Wentworth, {The asymptotics of the Arakelov-Green's function 
and Faltings' delta invariant}, {\em Comm. Math. Phys.}, {\bf 137} (1991)427-459 


\bibitem {Wp03} S. Wolpert, {Geometry of the {\WP} completion of 
{\TS}}, {\em Surveys in Differential Geometry, VIII: Papers in Honor 
of Calabi, Lawson, Siu and Uhlenbeck}, {Inter. Press} 2003

\bibitem {Ya} A. Yamada, {Precise variational formulas for {\Ad}s}, 
{\em Kodai. Math. J.} {\bf 3}, (1980) no. 1, 114-143

\end{thebibliography}
\end{document}